\documentclass[12pt]{article}
\usepackage{amsfonts}
\usepackage{amsmath}
\usepackage{amssymb}
\usepackage{graphicx}

\begin{document}

\begin{center}
{\bf Solution to Mathematics Magazine Problem 2065} \\
Yawen Zhang \\
Elizabethtown College \\
One Alpha Drive \\
Elizabethtown, PA 17022 \\
zhangy@etown.edu \\
\end{center}

\thispagestyle{empty}

\noindent
{\bf 2065.}
Let $\mathcal{Q}$ be a cube centered at the origin of $\mathbb{R}^3$.  Choose a unit vector $(a,b,c)$ uniformly at random on the surface of the unit sphere $a^2+b^2+c^2=1$, and let $\Pi$ be the plane $ax+by+cz=0$ through the origin and normal to $(a,b,c)$.  What is the probability that the intersection of $\Pi$ with $\mathcal{Q}$ is a hexagon?\cite{doi:10.1080/0025570X.2019.1544816}

\ \\

\noindent
{\bf Solution.} 
Without loss of generality, assume that $\mathcal{Q}$ has sides of length $2$.  Since $\mathcal{Q}$ has six faces, the intersection of $\mathcal{Q}$ and $\Pi$ is a hexagon only if $\Pi$ intersects two distinct edges on each of the six different faces of $\mathcal{Q}$. There are eight symmetric octants, so we focus on the first octant and then multiply by $8$ to account for the other octants.  If $a,b,c>0$, then the six edge intersections of $\Pi$ and $\mathcal{Q}$ are
$$
\pm \left (\frac {b-c}a, -1, 1 \right ), \quad 
\pm \left (-1, \frac {a-c}b, 1 \right ), \quad
\mbox{and} \quad
\pm \left (-1, 1, \frac {a-b}c \right ),
$$
where
$a<b+c$, $b<a+c$, $c<a+b$, and $a^2+b^2+c^2=1$.

\begin{center}
\includegraphics[height=3in]{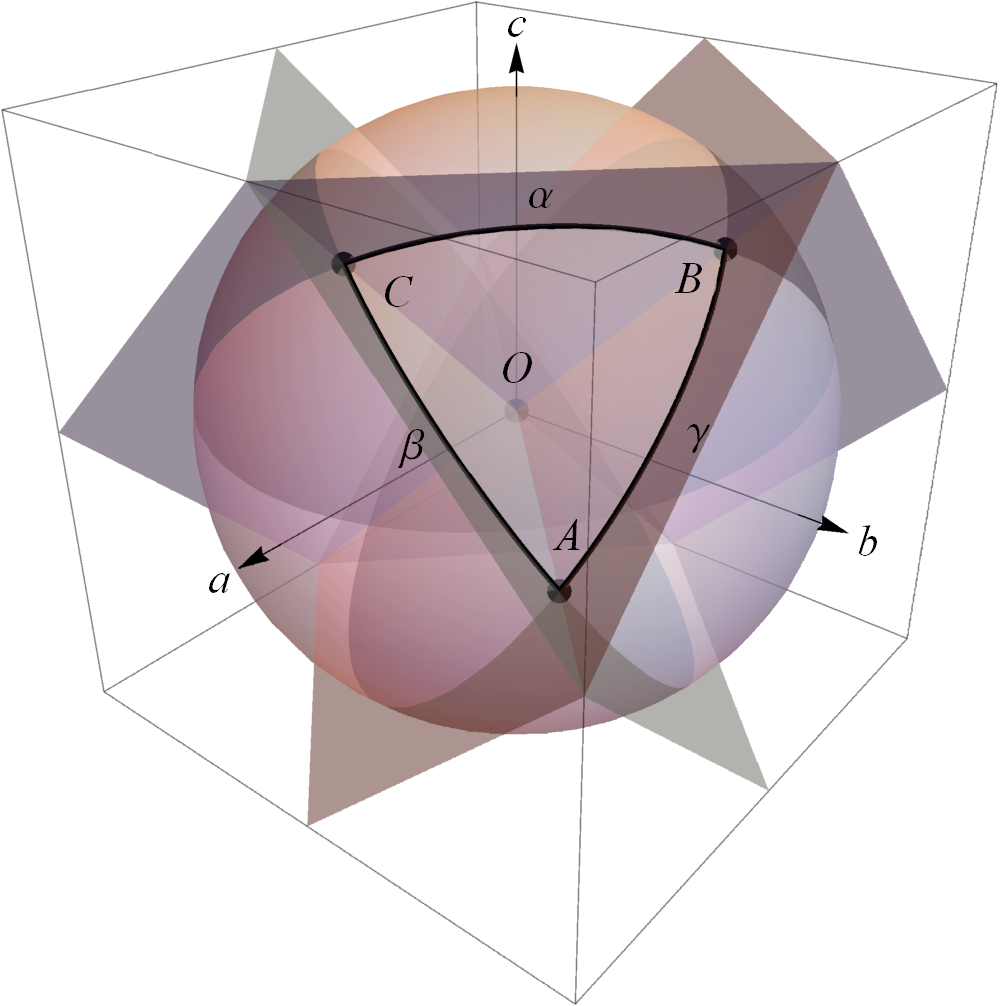}
\end{center}

\noindent
The planes $a=b+c$, $b=a+c$, and $c=a+b$ intersect $a^2+b^2+c^2=1$ in great circles that intersect in the first octant at points
$$
A \left ( \frac 1{\sqrt{2}}, \frac 1{\sqrt{2}}, 0 \right ), \quad
B \left ( 0, \frac 1{\sqrt{2}}, \frac 1{\sqrt{2}} \right ), \quad
\mbox{and} \quad
C \left ( \frac 1{\sqrt{2}}, 0, \frac 1{\sqrt{2}} \right ),
$$
as shown in the figure.  Choices of $(a,b,c)$ in the first octant that lead to hexagonal intersections of $\Pi$ and $\mathcal{Q}$ are in the spherical triangle $\Delta ABC$.

Let $O$ be the point at the origin and let $\alpha$, $\beta$, and $\gamma$ be the lengths of the arcs of $\Delta ABC$ opposite the angles at points $A$, $B$, and $C$, respectively.  Since $\overrightarrow{OA} \cdot \overrightarrow{OB} = \overrightarrow{OA} \cdot \overrightarrow{OC} = \overrightarrow{OB} \cdot \overrightarrow{OC} =1/2$ and the radius of the sphere is $1$, $\measuredangle{AOB} = \measuredangle{AOC} = \measuredangle{BOC} = \pi/3$ and $\alpha=\beta=\gamma=\pi/3$.  The spherical law of cosines states, for example, that 
$$
\cos \alpha = \cos \beta \cos \gamma + \sin \beta \sin \gamma \cos \left( \sphericalangle BAC\right ),
$$
so $\cos \left ( \sphericalangle BAC \right ) = \cos \left ( \sphericalangle ABC \right ) = \cos \left ( \sphericalangle ACB \right ) = 1/3$. The radius of the sphere is $1$, so the resulting area of $\Delta ABC$ is the sum of the interior angles minus $\pi$, or
$$
3\cos^{-1}  \left ( \frac {1}{3} \right ) - \pi.
$$
Since there are spherical $8$ triangles like $\Delta ABC$ and the surface area of the unit sphere is $4\pi$, the probability that the intersection of $\Pi$ and $\mathcal{Q}$ is a hexagon is
$$
\frac{8}{4\pi} \left[ 3\cos^{-1}  \left ( \frac {1}{3} \right ) - \pi \right ] = \frac 6{\pi} \cos^{-1} \left ( \frac 13 \right ) -2 \approx 0.351.
$$

\bibliographystyle{plain}
\bibliography{ref}

\end{document}